\documentclass{elsart}
\usepackage{amsfonts}
\usepackage{mathrsfs}
\usepackage{amssymb}
\date{}
\begin{document}
\begin{frontmatter}
\title{  Nonlinear maps preserving the mixed Jordan triple $\eta$-$*$-product between factors}
\thanks{This research was supported by the National Natural
  Science Foundation of China (No. 11601420) and the Natural
  Science Basic Research Plan in Shaanxi Province of China (Program No. 2018JM1053).\\        \indent $^*$ Corresponding
  author.
}
\author{Fangjuan Zhang}
\ead{zhfj888@126.com}
\address{School of Science, Xi'an University of Posts and Telecommunications, Xi'an 710121, P. R China }

\begin{abstract}   Let  $\mathcal{A}$ and $\mathcal{B}$ be two  factor von Neumann algebras and  $\eta$ be a  non-zero complex number.
A nonlinear bijective map $\phi:\mathcal
A\rightarrow\mathcal B$ has been demonstrated to satisfy $$\phi([A,B]_{*}^{\eta}\diamond_{\eta} C)=[\phi(A),\phi(B)]_{*}^{\eta}\diamond_{\eta}\phi(C)$$
 for all
$A,B,C\in\mathcal A.$ If $\eta=1,$ then $\phi$ is a linear $*$-isomorphism,  a conjugate linear $*$-isomorphism, the negative of a linear $*$-isomorphism,
or the negative of a conjugate linear $*$-isomorphism.
If $\eta\neq 1$ and satisfies $\phi(I)=1,$   then $\phi$ is  either a linear $*$-isomorphism or a conjugate linear $*$-isomorphism.

\vskip18pt \noindent{\it 2000AMS classification}  47B48; 46L10
\end{abstract}
\vskip -0.1in

\begin{keyword}  Mixed Jordan triple $\eta$-$*$-product ; Isomorphism;  von Neumann algebras

\end{keyword}
\end{frontmatter}

\section{Introduction}
 Let $\mathcal{A}$ be  a $*$-algebra. For a non-zero scalar $\eta,$  the Jordan $\eta$-$*$-product of two elements  $A,B\in\mathcal{A}$ is defined by
 $A\diamond_{\eta}B=AB+\eta BA^{*}.$
 The Jordan $(-1)$-$*$-product, customarily called the skew Lie product, has been extensively studied because it is naturally seen in the problem of representing quadratic functionals with sesquilinear functionals (\cite{12,13,14}) and that of characterizing ideals (\cite{10,11}). We often write the Jordan $(-1)$-$*$-product by $[A,B]_{*},$  the Jordan $(1)$-$*$-product by $A\bullet B$ and
  the Jordan $(-\eta)$-$*$-product by $[A,B]_{*}^{\eta},$ i.e., $[A,B]_{*}=AB-BA^{*}, A\bullet B=AB+BA^{*}, [A,B]_{*}^{\eta}=AB-\eta BA^{*}.$ A not necessarily linear map $\phi$
 between $*$-algebras $\mathcal{A}$ and $\mathcal{B}$ is said to preserve the Jordan $\eta$-$*$-product if $\phi(A\diamond_{\eta}B)=\phi(A)\diamond_{\eta}\phi(B)$
 for all $A,B\in\mathcal{A}.$ Recently, researchers have focused on  maps preserving the Jordan  $\eta$-$*$-product between $*$-algebra
  (\cite{1,2,3,8,9}).

  Huo et al. \cite{1} reported a more general problem, considering the Jordan triple $\eta$-$*$-product of three elements $A,B$ and $C$
 in a $*$-algebra $\mathcal{A}$ as defined by $A\diamond_{\eta}B\diamond_{\eta}C=(A\diamond_{\eta}B)\diamond_{\eta}C$ (note that $\diamond_{\eta}$ is not necessarily associative). Moreover, a map $\phi$ between $*$-algebras $\mathcal{A}$ and $\mathcal{B}$ is said to preserve the Jordan triple $\eta$-$*$-product if
 $\phi(A\diamond_{\eta}B\diamond_{\eta}C)=\phi(A)\diamond_{\eta}\phi(B)\diamond_{\eta}\phi(C)$ for all $A,B,C\in\mathcal{A}.$
 In \cite{1}, let $\eta\neq -1$ be a non-zero complex number and let $\phi$ be a bijection between two von Neumann algebras, one of which has no central abelian projections, satisfying
 $\phi(I)=I$ and preserving the Jordan triple $\eta$-$*$-product. Huo et al. reported that $\phi$ is a linear $*$-isomorphism if $\eta$ is not real and $\phi$ is the sum of a linear $*$-isomorphism and a conjugate linear $*$-isomorphism if $\eta$ is real. Nevertheless, Huo et al. have not considered the case $\eta=-1.$ Li et al. in  \cite{2} considered maps that preserve
 the Jordan triple $(-1)$-$*$-product without the assumption  $\phi(I)=I$  and confirmed that such a map between factors is  a linear $*$-isomorphism,  a conjugate linear $*$-isomorphism,  the negative of a linear $*$-isomorphism, or the negative of a conjugate linear $*$-isomorphism.  In \cite{3}, Zhao and Li discussed maps that preserve
 the Jordan triple $(1)$-$*$-product without the assumption  $\phi(I)=I$  and the results are similar to  \cite{2}.
The Lie product of two elements $A,B\in\mathcal{A}$ is defined by $[A,B]=AB-BA.$  A map $\phi$ between  factor von Neumann algebras   $\mathcal{A}$ and $\mathcal{B}$ is said to preserve Lie product if $\phi([A,B])=[\phi(A),\phi(B)]$ for all $A,B\in\mathcal{A}.$  In \cite{7}, Zhang and Zhang studied nonlinear bijective maps preserving Lie products
between factors.

The  mixed Jordan triple $\eta$-$*$-product  of three elements $A,B$ and $C$ in a $*$-algebra $\mathcal{A}$
defined by $[A,B]_{*}^{\eta}\diamond_{\eta}C.$ Moreover, a not necessarily linear map $\phi$ between factor von Neumann algebras  $\mathcal{A}$ and $\mathcal{B}$ preserves  mixed Jordan triple $\eta$-$*$-products
if $$\phi([A,B]_{*}^{\eta}\diamond_{\eta}C)=[\phi(A),\phi(B)]_{*}^{\eta}\diamond_{\eta}\phi(C)$$  for every $A,B,C\in\mathcal{A}.$ The primary result of this study confirms the following: Let
 $\mathcal{A}$ and $\mathcal{B}$ be two  factor von Neumann algebras and Let $\eta$ be a  non-zero complex number and $\phi:\mathcal{A}\rightarrow\mathcal{B}$ be a bijection that
 preserves the mixed Jordan triple $\eta$-$*$-product. Then  the following statements hold:\\
$($1$)$ If $\eta=1,$ then $\phi$ is a linear $*$-isomorphism,  a conjugate linear $*$-isomorphism,  the negative of a linear $*$-isomorphism, or the negative of a conjugate linear $*$-isomorphism.\\
$($2$)$ If $\eta\neq 1$  and satisfies $\phi(I)=1,$  then $\phi$ is  either a linear $*$-isomorphism or a conjugate linear $*$-isomorphism.

As usual, $\mathbb{R}$ and $\mathbb{C}$ denote respectively the real field and complex field. Throughout, algebras and spaces are over $\mathbb{C}.$
A von Neumann algebra $\mathcal{A}$ is a weakly closed, self-adjoint algebra of operators on a Hilbert space $H$ containing the identity operator $I.$
 $\mathcal{A}$ is a factor means that its center  only  contains the scalar operators. It is well known that the factor  $\mathcal{A}$ is prime, in the sense that $A\mathcal{A}B=\{0\}$ for $A,B\in\mathcal{A}$ implies either $A=0$ or $B=0.$
Let $\mathcal{P(A)}$ be the space of all projection operators of $\mathcal{A}.$

{\bf Lemma 1.1} ([3, Lemma 2.2])   {\it Let $\mathcal{A}$ be a  factor and $A\in\mathcal{A}.$ Then $AB+BA^{*}=0$ for all $B\in\mathcal{A}$ implies that $A\in\mbox{i}\mathbb{R}I$
($\mbox{i}$ is the imaginary number unit).}

{\bf Lemma 1.2} ([4, Lemma 2.2])   {\it Let $\mathcal{A}$ be a  factor von Neumann algebra and $A\in\mathcal{A}.$ If $[A,B]_{*}\in\mathbb{C}I$ for all $B\in\mathcal{A},$ then $A\in{\mathbb C}I.$}

{\bf Lemma 1.3} ([5, Problem 230])   {\it Let $\mathcal{A}$ be a Banach algebra with the identity  $I.$ If $A,B\in\mathcal{A}$ and $\lambda\in\mathbb{C}$ are such that $[A,B]=\lambda I,$ where $[A,B]=AB-BA,$ then $\lambda=0.$}

\section{ The mixed Jordan triple $1$-$*$-product preserving maps}

{\bf Theorem 2.1}   {\it Let $\mathcal{A},\mathcal{B}$ be two  factor von Neumann algebras. Suppose that $\phi$ is a bijective map from
 $\mathcal{A}$ to $\mathcal{B}$ with $\phi([A,B]_{*}\bullet C)=[\phi(A),\phi(B)]_{*}\bullet\phi(C)$ for all $A,B,C\in\mathcal A,$ then $\phi$  is a linear $*$-isomorphism, a conjugate linear $*$-isomorphism, the negative of a linear $*$-isomorphism, or the negative conjugate linear $*$-isomorphism.}

 First we give a key technique. Assume that $A_{1},A_{2},\cdots,A_{n}$ and $T$ are in $\mathcal{A}$ with $\phi(T)=\Sigma_{i=1}^{n}\phi(A_{i}).$
 Then for $S_{1},S_{2},S_{3}\in\mathcal{A},$ we obtain
 $$\phi([T,S_{2}]_{*}^{\eta}\diamond_{\eta}S_{3})=[\phi(T),\phi(S_{2})]_{*}^{\eta}\diamond_{\eta}\phi(S_{3})=\displaystyle\sum_{i = 1}^{n}\phi([A_{i},S_{2}]_{*}^{\eta}\diamond_{\eta}S_{3}),\eqno(1)$$
$$\phi([S_{1},T]_{*}^{\eta}\diamond_{\eta}S_{3})=[\phi(S_{1}),\phi(T)]_{*}^{\eta}\diamond_{\eta}\phi(S_{3})=\displaystyle\sum_{i = 1}^{n}\phi([S_{1},A_{i}]_{*}^{\eta}\diamond_{\eta} S_{3})\eqno(2)$$
and
$$\phi([S_{1},S_{2}]_{*}^{\eta}\diamond_{\eta}T)=[\phi(S_{1}),\phi(S_{2})]_{*}^{\eta}\diamond_{\eta}\phi(T)=\displaystyle\sum_{i = 1}^{n}\phi([S_{1},S_{2}]_{*}^{\eta}\diamond_{\eta} A_{i}).\eqno(3)$$

Choose an arbitrary nontrivial projection $P_{1}\in\mathcal{A},$ write $P_{2}=I-P_{1}.$ Denote $\mathcal{A}_{ij}=P_{i}\mathcal{A}P_{j},i,j,=1,2,$ then
 $\mathcal{A}=\sum_{i,i=1}^{2}\mathcal{A}_{ij}.$ For every $A\in\mathcal{A},$ we can write it as $A=\sum_{i,i=1}^{2}A_{ij},$ where $A_{ij}$ denotes an arbitrary element of
 $\mathcal{A}_{ij}.$ We will complete the proof of Theorem 2.1 by proving several lemmas.

 {\bf Lemma 2.1}  {\it $\phi(0)=0.$}

 {\bf Proof} Since $\phi$ is surjective, there exists $A\in\mathcal{A}$ such that $\phi(A)=0.$ Hence $\phi(0)=\phi([0,A]_{*}\bullet A)=[\phi(0),\phi(A)]_{*}\bullet\phi(A)=0.$

{\bf Lemma 2.2}  {\it $\phi(\mathbb{R}I)=\mathbb{R}I, \phi(\mathbb{C}I)=\mathbb{C}I$  and $\phi$ preserves self-adjoint elements in both directions.}

 {\bf Proof} Let $\lambda\in\mathbb{R}$ be arbitrary. It is easily seen that $$0=\phi([\lambda I,B]_{*}\bullet C)=[\phi(\lambda I),\phi(B)]_{*}\bullet\phi(C)$$ holds true for any $B,C\in\mathcal{A}.$
Since $\phi$ is surjective, by Lemma 1.1, which indicates that $$[\phi(\lambda I),\phi(B)]_{*}\in\mbox{i}\mathbb{R}I.$$  Then $[\phi(\lambda I),B]_{*}\in\mathbb{C}I$ for any $B\in\mathcal{B}.$
We obtain from Lemma 1.2 that $\phi(\lambda I)\in\mathbb{C}I,$ so there exists $\lambda_{0}\in{\mathbb{C}}$ such that $(\lambda_{0}-\overline{\lambda_{0}})B\in{\mathbb{C}}I$
for any  $B\in\mathcal{B},$ then $\phi(\lambda I)\in\mathbb{R}I.$ Note that $\phi^{-1}$ has the same properties as $\phi.$ Similarly, if $\phi(A)\in\mathbb{R}I,$
then $A\in\mathbb{R}I.$ Therefore, $\phi(\mathbb{R}I)=\mathbb{R}I.$

 Since $\phi(\mathbb{R}I)=\mathbb{R}I,$ exists $\lambda\in\mathbb{R}$ such that $\phi(\lambda I)=I.$ For any $A=A^{*}\in\mathcal{A}$ and $B\in\mathcal{A},$ we obtain
$$0=\phi([A,\lambda I]_{*}\bullet B)=[\phi(A),I]_{*}\bullet\phi(B),$$
by the surjectivity of  $\phi$ and Lemma 1.1, the above equation indicates $[\phi(A),I]_{*}\in\mbox{i}\mathbb{R}I.$ Then there exists $\lambda\in\mbox{i}\mathbb{R}$ such that $\phi(A)^{*}=\phi(A)+\lambda I.$ However, $$0=\phi([A,A]_{*}\bullet B)=[\phi(A),\phi(A)]_{*}\bullet\phi(B)$$ for all $A=A^{*}\in\mathcal{A}$ and $B\in\mathcal{A}.$
Similarly, $[\phi(A),\phi(A)]_{*}\in\mbox{i}\mathbb{R}I.$ Then $\lambda \phi(A)\in\mbox{i}\mathbb{R}I.$
If $\lambda\neq 0,$ then  $\phi(A)\in\mathbb{R}I.$ It follows from $\phi(\mathbb{R}I)=\mathbb{R}I$ that $A=A^{*}\in\mathbb{R}I,$
which is contradiction. Thus $\lambda=0.$ Now we get that $\phi(A)=\phi(A)^{*}.$
Similarly, if $\phi(A)=\phi(A)^{*},$ then $A=A^{*}\in\mathcal{A}.$ Therefore $\phi$ preserves self-adjoint elements in both directions.

Let $\lambda\in\mathbb{C}$ be arbitrary. For every  $A=A^{*}\in\mathcal{A},$ we obtain that
$$0=\phi([A,\lambda I]_{*}\bullet B)=[\phi(A),\phi(\lambda I)]_{*}\bullet\phi(B)$$
holds true for any  $B\in\mathcal{A}.$ By the surjectivity of $\phi$ and Lemma 1.1 again, the above equation indicates $[\phi(A),\phi(\lambda I)]_{*}\in\mbox{i}\mathbb{R}I.$
Since $A=A^{*},$ we have $\phi(A)=\phi(A)^{*}.$ Hence $[\phi(A),\phi(\lambda I)]\in\mbox{i}\mathbb{R}I.$
We obtain from Lemma 1.3 that  $[\phi(A),\phi(\lambda I)]=0,$ and then $B\phi(\lambda I)=\phi(\lambda I)B$ for any $B=B^{*}\in\mathcal{B}.$
Thus for any  $B\in\mathcal{B},$ since $B=B_{1}+\mbox{i}B_{2}$ with $B_{1}=\frac{B+B^{*}}{2}$ and  $B_{2}=\frac{B-B^{*}}{2\mbox{i}},$ we get
 $$B\phi(\lambda I)=\phi(\lambda I)B$$ for any $B\in\mathcal{B}.$ Hence $\phi(\lambda I)\in\mathbb{C}I.$
Similarly, if $\phi(A)\in\mathbb{C}I,$ then $A\in\mathbb{C}I.$ Therefore, $\phi(\mathbb{C}I)=\mathbb{C}I.$

 {\bf Lemma 2.3}  {\it $\phi(\mathcal{P(A)}+\mathbb{R}I)=\mathcal{P(B)}+\mathbb{R}I.$ }

 {\bf Proof} Fix a nontrivial projection $P\in\mathcal{P(B)}.$ Based on Lemma 2.2, exists  $A=A^{*}\in\mathcal{A}$ such that $\phi(A)=P+\mathbb{R}I.$
For any $B=B^{*}\in\mathcal{A}$ and $C\in\mathcal{A},$ we then have the following:
\begin{eqnarray*}
&&\phi([A,B]_{*}\bullet C)=[\phi(A),\phi(B)]_{*}\bullet\phi(C)\\
&=&[P,\phi(B)]_{*}\bullet\phi(C)=[([P,\phi(B)]_{*}\bullet P),P]_{*}\bullet\phi(C)\\
&=&[([\phi(A),\phi(B)]_{*}\bullet\phi(A)),\phi(A)]_{*}\bullet\phi(C)=\phi([([A,B]_{*}\bullet A),A]_{*}\bullet C).
\end{eqnarray*}
By the injectivity of $\phi,$ we concur that $[([A,B]_{*}\bullet A),A]_{*}\bullet C=[A,B]_{*}\bullet C,$ which indicates $[([A,B]_{*}\bullet A),A]_{*}-[A,B]_{*}\in\mbox{i}\mathbb{R}I$ for all $B=B^{*}\in\mathcal{A}.$
For every $X\in\mathcal{A},$ we have $X=X_{1}+\mbox{i}X_{2},$ where $X_{1}=\frac{X+X^{*}}{2}$ and  $X_{2}=\frac{X-X^{*}}{2\mbox{i}}$ are self-adjoint.
We obtain $[A,[A,[A,X]]]-[A,X]\in\mbox{i}\mathbb{R}I,$ i.e.,
$$A^{3}X-3A^{2}XA+3AXA^{2}-XA^{3}-AX+XA\in\mbox{i}\mathbb{R}I\eqno(4)$$
holds true for any  $X\in\mathcal{A}.$

Let $\mathcal{U}$ be the group of unitary operators of $\mathcal{A}$ and let $\varphi$ be the set of the functions $U\rightarrow f(U)$ defined on $\mathcal{U}$
with non-negative real values, zero except on a finite subset of  $\mathcal{U}$ and such that $\sum_{U\in\mathcal{U}}f(U)=1.$ For $A\in\mathcal{A}$ and $f\in\varphi,$
we put $f\cdot A=\sum_{U\in\mathcal{U}}f(U)UAU^{*}.$

For any $U\in\mathcal{U},$  by Eq. (4),
$$(A^{3}-A)U-3A^{2}UA+3AUA^{2}-U(A^{3}-A)=\alpha I$$
for certain $\alpha\in\mbox{i}\mathbb{R}I.$ It follows that
$$A^{3}-A-3A^{2}UAU^{*}+3AUA^{2}U^{*}-U(A^{3}-A)U^{*}=\alpha U^{*},$$
and so $A^{3}-A-3A^{2}f\cdot A+3Af\cdot A^{2}-f\cdot A^{3}+f\cdot A=\alpha U^{*}$ for any $f\in\varphi.$ Since $\mathcal{A}$ is a factor von Neumann algebra, we obtain from
 [6, Lemma 5 (Part $\mbox{\uppercase\expandafter{\romannumeral3}},$ Chapter 5)] that there exist $\lambda_{1},\lambda_{2},\lambda_{3}\in\mathbb{C}$ such that
 $$A^{3}-A-3\lambda_{1} A^{2}+3\lambda_{2} A-(\lambda_{3}-\lambda_{1})I=\alpha U^{*}.$$
Thus $U(A^{3}-A)U^{*}-3\lambda_{1} UA^{2}U^{*}+3\lambda_{2} UAU^{*}-(\lambda_{3}-\lambda_{1})I=\alpha U^{*}$
and then $f\cdot A^{3}-f\cdot A-3\lambda_{1} f\cdot A^{2}+3\lambda_{2} f\cdot A-(\lambda_{3}-\lambda_{1})I=\alpha U^{*}$ for any $f\in\varphi.$
By [6, Lemma 5 (Part $\mbox{\uppercase\expandafter{\romannumeral3}},$ Chapter 5)] again, we have $\alpha U^{*}=0$ for any $U\in\mathcal{U}.$
Hence $\alpha=0.$ Thus we have
$$(A^{3}-A)U-3A^{2}UA+3AUA^{2}-U(A^{3}-A)=0\eqno(5)$$
and
$$A^{3}-A=3\lambda_{1} A^{2}-3\lambda_{2} A+(\lambda_{3}-\lambda_{1})I\eqno(6)$$
 for any  $U\in\mathcal{U}.$ By Eqs. (5)--(6), we conclude that
$$(\lambda_{1}A^{2}-\lambda_{2}A)U-A^{2}UA+AUA^{2}-U(\lambda_{1}A^{2}-\lambda_{2}A)=0$$
and
$$(\lambda_{1}A^{2}-\lambda_{2}A)UAU^{*}-A^{2}UA^{2}U^{*}+AUA^{3}U^{*}-U(\lambda_{1}A^{2}-\lambda_{2}A)AU^{*}=0$$
 for any $U\in\mathcal{U}.$ Thus
$$(\lambda_{1}A^{2}-\lambda_{2}A)f\cdot A-A^{2}f\cdot A^{2}+Af\cdot A^{3}-\lambda_{1}f\cdot A^{3}+\lambda_{2}f\cdot A^{2}=0$$
 for any $f\in\varphi.$ By applying [6, Lemma 5 (Part $\mbox{\uppercase\expandafter{\romannumeral3}},$ Chapter 5)] again, we obtain
 $$\lambda_{1}(\lambda_{1}A^{2}-\lambda_{2}A)-\lambda_{2}A^{2}+\lambda_{3}A+(\lambda_{2}^{2}-\lambda_{1}\lambda_{3})I=0,$$
i.e.,
$$(\lambda_{1}^{2}-\lambda_{2})A^{2}+(\lambda_{3}-\lambda_{1}\lambda_{2})A+(\lambda_{2}^{2}-\lambda_{1}\lambda_{3})I=0.\eqno(7)$$

 If $\lambda_{2}=\lambda_{1}^{2},$ we obtain from $\phi(\mathbb{C}I)=\mathbb{C}I$ and $\phi(A)=P+\mathbb{R}I\notin\mathbb{C}I$ that $A\notin\mathbb{C}I,$ then
 $\lambda_{3}=\lambda_{1}\lambda_{2}=\lambda_{1}^{3},$ so by Eq. (6), we have $(A-\lambda_{1} I)^{3}=A-\lambda_{1} I.$ Let $B=A-\lambda_{1} I,$ then
$$B^{3}=B \ \ \mbox{and} \ \ [B,[B,[B,X]]]=[B,X]$$
for any $X\in\mathcal{A}.$ This indicates that
$$B^{2}XB-BXB^{2}=0\eqno(8)$$
for any $X\in\mathcal{A}.$ Let $E_{1}=\frac{1}{2}(B^{2}+B)$ and $E_{2}=\frac{1}{2}(B^{2}-B).$ It follows from the fact $B^{3}=B$ that $E_{1}$ and $E_{2}$
are idempotents of $\mathcal{A}$ and that
$$B=E_{1}-E_{2}, \ \ B^{2}=E_{1}+E_{2}, \ \ E_{1}E_{2}=E_{2}E_{1}=0.$$
This along with Eq. (8) shows us that $E_{1}XE_{2}=0$ for any $X\in\mathcal{A}.$ Then $E_{1}=0$ or $E_{2}=0,$ and so $A=\lambda_{1}I+E_{1}$ or $A=\lambda_{1}I-E_{2}$
is the sum of a scalar and an idempotent of $\mathcal{A}.$

 If $\lambda_{2}\neq\lambda_{1}^{2},$ by Eq. (7), we have $A^{2}=\lambda A+\mu I$ for certain $\lambda,\mu\in\mathbb{C}.$ This along with Eq. (5) indicates that
$$(\lambda^{2}+4\mu-1)(AU-UA)=0\eqno(9)$$
for all $U\in\mathcal{U}.$ Since $A\notin\mathbb{C}I,$ we have $AU-UA\neq 0$ for some $U\in\mathcal{U}.$
By Eq. (9), we obtain $\lambda^{2}+4\mu-1=0.$ Let $E=A+\frac{1}{2}(1-\lambda)I,$ then
\begin{eqnarray*}
E^{2}&=&A^{2}+(1-\lambda)A+\frac{1}{4}(1-\lambda)^{2}I=\lambda A+\mu I+(1-\lambda)A+\frac{1}{4}(1-\lambda)^{2}I\\
&=&A+\frac{1}{4}(\lambda^{2}+4\mu-2\lambda+1)I=A+\frac{1}{2}(1-\lambda)I=E.
\end{eqnarray*}
Hence $A=\frac{1}{2}(\lambda-1)I+E$ is the sum of a scalar and an idempotent of $\mathcal{A}.$ Since $A=A^{*},$ then  $A=\alpha I+E, \alpha\in\mathbb{R}I, E\in\mathcal{P(A)}.$
If $E=0$ or  $E=I,$ from $ \phi(A)=P+\mathbb{R}I,$ we obtain $\phi(\mathbb{R}I)=P+\mathbb{R}I.$ It follows $\phi(\mathbb{R}I)=\mathbb{R}I$ that $P=0$ or $P=I,$ since $P$
is a nontrivial projection, which is a contradiction.
Thus, $A$  is the sum of a real number and  a nontrivial projection  of $\mathcal{A}.$ Applying the same argument to $\phi^{-1},$ we can obtain the reverse inclusion and equality follows.

{\bf Lemma 2.4} {\it  Let $i,j\in\{1,2\}$ with $i\neq j.$ Then $\phi(A_{ii}+B_{ij})=\phi(A_{ii})+\phi(B_{ij})$ and $\phi(A_{ii}+B_{ji})=\phi(A_{ii})+\phi(B_{ji})$
 for all $A_{ii}\in\mathcal{A}_{ii}, B_{ij}\in\mathcal{A}_{ij}$ and $B_{ji}\in\mathcal{A}_{ji}.$}

 {\bf Proof}  Let $X=\sum_{i,i=1}^{2}X_{ij}\in\mathcal{A}$ such that $\phi(X)=\phi(A_{ii})+\phi(B_{ij}).$ It follows from Eq. (2) that
$$\phi(X_{ij}+X_{ij}^{*})=\phi([P_{i},X]_{*}\bullet P_{j})=\phi([P_{i},A_{ii}]_{*}\bullet P_{j})+\phi([P_{i},B_{ij}]_{*}\bullet P_{j})=\phi(B_{ij}+B_{ij}^{*}).$$
By the injectivity of $\phi,$ we have $X_{ij}=B_{ij}.$
We obtain from Eq. (2) that
$$\phi(X_{ji}+X_{ji}^{*})=\phi([P_{j},X]_{*}\bullet P_{i})=\phi([P_{j},A_{ii}]_{*}\bullet P_{i})+\phi([P_{j},B_{ij}]_{*}\bullet P_{i})=0,$$
which indicates that
$X_{ji}=0.$
For every $T_{ij}\in\mathcal{A}_{ij},$ by applying Eq. (2) again, we obtain
\begin{eqnarray*}
&&\phi(T_{ij}X_{jj}+X_{jj}^{*}T_{ij}^{*})=\phi([T_{ij},X]_{*}\bullet P_{j})\\
&=&\phi([T_{ij},A_{ii}]_{*}\bullet P_{j})+\phi([T_{ij},B_{ij}]_{*}\bullet P_{j})=0,
\end{eqnarray*}
which implies that $T_{ij}X_{jj}=X_{jj}^{*}T_{ij}^{*}=0$ for all $T_{ij}\in \mathcal{A}_{ij}.$ By the primeness of $\mathcal{A},$
we get  $X_{jj}=0.$
For every $T_{ij}\in\mathcal{A}_{ij},$  by applying Eq. (1), we obtain
\begin{eqnarray*}
&&\phi(X_{ii}T_{ij}+T_{ij}^{*}X_{ii}^{*})=\phi([X,T_{ij}]_{*}\bullet P_{j})\\
&=&\phi([A_{ii},T_{ij}]_{*}\bullet P_{j})+\phi([B_{ij},T_{ij}]_{*}\bullet P_{j})=\phi(A_{ii}T_{ij}+T_{ij}^{*}A_{ii}^{*}),
\end{eqnarray*}
which indicates that $X_{ii}T_{ij}=A_{ii}T_{ij}$ for all $T_{ij}\in \mathcal{A}_{ij}.$ By the primeness of $\mathcal{A},$
we obtain $X_{ii}=A_{ii}.$ Thus $\phi(A_{ii}+B_{ij})=\phi(A_{ii})+\phi(B_{ij}).$
In the second case, we can similarly prove that the conclusion is valid.

{\bf Lemma 2.5} {\it  Let $i,j\in\{1,2\}$ with $i\neq j.$ Then $\phi(A_{ij}+B_{ji})=\phi(A_{ij})+\phi(B_{ji})$  for all $A_{ij}\in\mathcal{A}_{ij}$ and $B_{ji}\in\mathcal{A}_{ji}.$}

 {\bf Proof}   Choose $X=\sum_{i,i=1}^{2}X_{ij}\in\mathcal{A}$ such that $\phi(X)=\phi(A_{ij})+\phi(B_{ji}).$  It follows from Eq. (1) that
\begin{eqnarray*}
&&\phi(X_{ji}+X_{ji}^{*})=\phi([X,P_{i}]_{*}\bullet P_{i})\\
&=&\phi([A_{ij},P_{i}]_{*}\bullet P_{i})+\phi([B_{ji},P_{i}]_{*}\bullet P_{i})=\phi(B_{ji}+B_{ji}^{*}),
\end{eqnarray*}
Thus we have $X_{ji}=B_{ji}.$
Similarly, $X_{ij}=A_{ij}.$  For every $T_{ij}\in\mathcal{A}_{ij},$  By applying Eq. (2), we obtain
 \begin{eqnarray*}
&&\phi(T_{ij}X_{jj}+X_{jj}^{*}T_{ij}^{*})=\phi([T_{ij},X]_{*}\bullet P_{j})\\
&=&\phi([T_{ij},A_{ij}]_{*}\bullet P_{j})+\phi([T_{ij},B_{ji}]_{*}\bullet P_{j})=0,
\end{eqnarray*}
from this, we get $X_{jj}=0.$ In the same manner, we obtain  $X_{ii}=0.$

{\bf Lemma 2.6} {\it  $\phi(\Sigma_{i,j=1}^{2}A_{ij})=\Sigma_{i,j=1}^{2}\phi(A_{ij})$  for all $A_{ij}\in \mathcal{A}_{ij}.$}

 {\bf Proof}   Let $X=\sum_{i,i=1}^{2}X_{ij}\in\mathcal{A}$ such that $\phi(X)=\Sigma_{i,j=1}^{2}\phi(A_{ij}).$
 It follows from Eq. (2) that
$\phi([P_{1},X]_{*}\bullet P_{2})=\Sigma_{i,j=1}^{2}\phi([P_{1},A_{ij}]_{*}\bullet P_{2}),$
i.e.,
$\phi(X_{12}+X_{12}^{*})=\phi(A_{12}+A_{12}^{*}),$
which implies that $X_{12}=A_{12}.$ Similarly,  $X_{21}=A_{21}.$
For every  $T_{12}\in \mathcal{A}_{12},$  by applying Eq. (2) again, we obtain
$\phi([T_{12},X]_{*}\bullet P_{2})=\Sigma_{i,j=1}^{2}\phi([T_{12},A_{ij}]_{*}\bullet P_{2})$
 for all $T_{12}\in \mathcal{A}_{12}.$
Thus we have $X_{22}=A_{22}.$
In the same manner, we obtain $X_{11}=A_{11}.$

{\bf Lemma 2.7} {\it  Let $i,j\in\{1,2\}$ with $i\neq j.$ Then $\phi(A_{ij}+B_{ij})=\phi(A_{ij})+\phi(B_{ij})$  for all $A_{ij}\in\mathcal{A}_{ij}$ and $B_{ij}\in\mathcal{A}_{ij}.$}

{\bf Proof}  It follows from $A_{ij}+B_{ij}+A_{ij}^{*}+B_{ij}A_{ij}^{*}=[-\frac{\mbox{i}}{2}I,\mbox{i}P_{i}+\mbox{i}A_{ij}]_{*}\bullet(P_{j}+B_{ij})$   and  Lemmas 2.6, 2.4, 2.5 that
\begin{eqnarray*}
&&\phi(A_{ij}+B_{ij})+\phi(A_{ij}^{*})+\phi(B_{ij}A_{ij}^{*})\\
&=&\phi(A_{ij}+B_{ij}+A_{ij}^{*}+B_{ij}A_{ij}^{*})\\
&=&\phi([-\frac{\mbox{i}}{2}I,\mbox{i}P_{i}+\mbox{i}A_{ij}]_{*}\bullet(P_{j}+B_{ij}))\\
&=&[\phi(-\frac{\mbox{i}}{2}I),\phi(\mbox{i}P_{i}+\mbox{i}A_{ij})]_{*}\bullet\phi(P_{j}+B_{ij})\\
&=&[\phi(-\frac{\mbox{i}}{2}I),\phi(\mbox{i}P_{i})+\phi(\mbox{i}A_{ij})]_{*}\bullet(\phi(P_{j})+\phi(B_{ij}))\\
&=&\phi([-\frac{\mbox{i}}{2}I,\mbox{i}P_{i}]_{*}\bullet P_{j})+\phi([-\frac{\mbox{i}}{2}I,\mbox{i}P_{i}]_{*}\bullet B_{ij})\\
&&+\phi([-\frac{\mbox{i}}{2}I,\mbox{i}A_{ij}]_{*}\bullet P_{j})+\phi([-\frac{\mbox{i}}{2}I,\mbox{i}A_{ij}]_{*}\bullet B_{ij})\\
&=&\phi(B_{ij})+\phi(A_{ij}+A_{ij}^{*})+\phi(B_{ij}A_{ij}^{*})\\
&=&\phi(B_{ij})+\phi(A_{ij})+\phi(A_{ij}^{*})+\phi(B_{ij}A_{ij}^{*}),
\end{eqnarray*}
which indicates that $\phi(A_{ij}+B_{ij})=\phi(A_{ij})+\phi(B_{ij}).$

{\bf Lemma 2.8} {\it  Let $i\in\{1,2\}.$  Then $\phi(A_{ii}+B_{ii})=\phi(A_{ii})+\phi(B_{ii})$  for all $A_{ii}\in\mathcal{A}_{ii}$ and $B_{ii}\in\mathcal{A}_{ii}.$}

 {\bf Proof}  Choose $X=\sum_{i,i=1}^{2}X_{ij}\in\mathcal{A}$  such that $\phi(X)=\phi(A_{ii})+\phi(B_{ii}).$  It follows from Eq. (2) that
\begin{eqnarray*}
&&\phi(X_{ij}+X_{ij}^{*})=\phi([P_{i},X]_{*}\bullet P_{j})\\
&=&\phi([P_{i},A_{ii}]_{*}\bullet P_{j})+\phi([P_{i},B_{ii}]_{*}\bullet P_{j})=0.
\end{eqnarray*}
Thus we have  $X_{ij}=0.$  Similarly, $X_{ji}=0.$
For every $T_{ij}\in\mathcal{A}_{ij},$  By applying Eq. (2) again, we have
\begin{eqnarray*}
&&\phi(T_{ij}X_{jj}+X_{jj}^{*}T_{ij}^{*})=\phi([T_{ij},X]_{*}\bullet P_{j})\\
&=&\phi([T_{ij},A_{ii}]_{*}\bullet P_{j})+\phi([T_{ij},B_{ii}]_{*}\bullet P_{j})=0,
\end{eqnarray*}
 which implies that $T_{ij}X_{jj}=X_{jj}^{*}T_{ij}^{*}=0.$  By the primeness of $\mathcal{A},$
we obtain $X_{jj}=0.$ Therefore,
$$\phi(X_{ii})=\phi(A_{ii})+\phi(B_{ii}).\eqno(10)$$
For every $T_{ij}\in\mathcal{A}_{ij},$ it follows from  Eq. (1) and Lemmas 2.5 and 2.7 that
\begin{eqnarray*}
&&\phi(X_{ii}T_{ij}+T_{ij}^{*}X_{ii}^{*})=\phi([X,T_{ij}]_{*}\bullet P_{j})\\
&=&\phi([A_{ii},T_{ij}]_{*}\bullet P_{j})+\phi([B_{ii},T_{ij}]_{*}\bullet P_{j})\\
&=&\phi(A_{ii}T_{ij}+T_{ij}^{*}A_{ii}^{*})+\phi(B_{ii}T_{ij}+T_{ij}^{*}B_{ii}^{*})\\
&=&\phi(A_{ii}T_{ij})+\phi(T_{ij}^{*}A_{ii}^{*})+\phi(B_{ii}T_{ij})+\phi(T_{ij}^{*}B_{ii}^{*})\\
&=&\phi(A_{ii}T_{ij}+B_{ii}T_{ij})+\phi(T_{ij}^{*}A_{ii}^{*}+T_{ij}^{*}B_{ii}^{*})\\
&=&\phi(A_{ii}T_{ij}+B_{ii}T_{ij}+T_{ij}^{*}A_{ii}^{*}+T_{ij}^{*}B_{ii}^{*}),
\end{eqnarray*}
which indicates that $X_{ii}=A_{ii}+B_{ii}.$
  This together with Eq. (10)  shows that $\phi(A_{ii}+B_{ii})=\phi(A_{ii})+\phi(B_{ii}).$

{\bf Lemma 2.9} {\it $\phi$ is additive and $\phi(\mathcal{P(A)})=\mathcal{P(B)}.$}

 {\bf Proof}   By Lemmas 2.6--2.8, $\phi$ is additive.
By Lemmas 2.3 and 2.2, we have $\phi(\mathcal{P(A)})=\mathcal{P(B)}.$

{\bf Remark 2.1}  Since $[P_{1},B]_{*}\bullet C=[B,P_{2}]_{*}\bullet C$ for all $B=B^{*}\in\mathcal{A}$ and $C\in{\mathcal{A}},$
 from Lemma 2.9  $$[Q_{1},\phi(B)]_{*}\bullet\phi(C)=[\phi(B),Q_{2}]_{*}\bullet\phi(C),$$ where  $Q_{i}\in\mathcal{P(B)}, i=1,2.$
The surjectivity of $\phi$ indicates that $[Q_{1},\phi(B)]_{*}-[\phi(B),Q_{2}]_{*}\in\mbox{i}\mathbb{R}I.$ It follows from Lemma 2.2 that
$[Q_{1}+Q_{2},B]\in\mbox{i}\mathbb{R}I$ holds true for all  $B=B^{*}\in\mathcal{B}.$  By Lemma 1.3, $[Q_{1}+Q_{2},B]=0.$
 Thus for every  $B\in\mathcal{B},$ because $B=B_{1}+\mbox{i}B_{2}$ with $B_{1}=\frac{B+B^{*}}{2}$ and  $B_{2}=\frac{B-B^{*}}{2\mbox{i}},$ we get
$[Q_{1}+Q_{2},B]=0$ for all $B\in\mathcal{B}.$ From this, exists $\lambda\in\mathbb{R}$ such that $$Q_{1}+Q_{2}=\lambda I.$$

Multiplying by $Q_{1}$ and $Q_{2}$ from the left and right respectively in the above equation, we obtain $Q_{1}+Q_{1}Q_{2}=\lambda Q_{1}$ and
$Q_{1}Q_{2}+Q_{2}=\lambda Q_{2}.$ Therefore, we can concur that $(1-\lambda)(Q_{1}-Q_{2})=0$ by subtracting the above two equations. By the injectivity of $\phi,$
 exists $P_{1}\neq P_{2}$ such that $Q_{1}\neq Q_{2}.$ Thus $\lambda=1$ and then $Q_{2}=I-Q_{1}.$

{\bf Lemma 2.10} {\it  $\phi(\mathcal{A}_{ij})=\mathcal{B}_{ij}, i,j=1,2.$}

 {\bf Proof}  Let $i,j\in\{1,2\}$ with $i\neq j$ and $A_{ij}\in\mathcal{A}_{ij}.$  By the fact $A_{ij}=[P_{i},A_{ij}]_{*}\bullet \frac{I}{2},$ Lemma 2.9 and  Remark 2.1, we have
$$\phi(A_{ij})=(Q_{i}\phi(A_{ij})-\phi(A_{ij})Q_{i})\phi(\frac{I}{2})-\phi(\frac{I}{2})(\phi(A_{ij})^{*}Q_{i}-Q_{i}\phi(A_{ij})^{*}).$$ From this and Lemma 2.2, we get
$Q_{i}\phi(A_{ij})Q_{i}=Q_{j}\phi(A_{ij})Q_{j}=0.$
Thus
$$\phi(A_{ij})=Q_{i}\phi(A_{ij})Q_{j}+Q_{j}\phi(A_{ij})Q_{i}.\eqno(11)$$
For every $B\in\mathcal{A},$ we obtain from the fact $[A_{ij},P_{i}]_{*}\bullet B=0,$ Lemma 2.9  and  Remark 2.1 that $[\phi(A_{ij}),Q_{i}]_{*}\bullet\phi(B)=0.$
Thus $[\phi(A_{ij}),Q_{i}]_{*}\in\mbox{i}\mathbb{R}I,$ which together with Eq. (11)  indicates that
$Q_{j}\phi(A_{ij})Q_{i}-Q_{i}\phi(A_{ij})^{*}Q_{j}\in\mbox{i}\mathbb{R}I.$ Multiplying by $Q_{j}$ and $Q_{i}$ from the left and right respectively in the above equation, we have
$Q_{j}\phi(A_{ij})Q_{i}=0.$ It follows from Eq. (11) that $\phi(A_{ij})=Q_{i}\phi(A_{ij})Q_{j},$ and then we obtain $\phi({\mathcal A}_{ij})\subseteq{\mathcal B}_{ij}.$
 Applying the same argument to $\phi^{-1},$ we obtain ${\mathcal B}_{ij}\subseteq\phi({\mathcal A}_{ij}).$ Thus $\phi({\mathcal A}_{ij})={\mathcal B}_{ij}, i\neq j.$

 Let $A_{jj}\in\mathcal{A}_{jj}$ and $i\neq j.$ It follows from  Lemma 2.9 and  Remark 2.1 that
$$0=\phi([P_{i},A_{jj}]_{*}\bullet P_{j})=[Q_{i},\phi(A_{jj})]_{*}\bullet Q_{j}=Q_{i}\phi(A_{jj})Q_{j}+Q_{j}\phi(A_{jj})^{*}Q_{i}$$
and
$$0=\phi([P_{j},A_{jj}]_{*}\bullet P_{i})=[Q_{j},\phi(A_{jj})]_{*}\bullet Q_{i}=Q_{j}\phi(A_{jj})Q_{i}+Q_{i}\phi(A_{jj})^{*}Q_{j}.$$
which indicates that $Q_{i}\phi(A_{jj})Q_{j}=Q_{j}\phi(A_{jj})Q_{i}=0.$ Now we obtain
$$\phi(A_{jj})=Q_{i}\phi(A_{jj})Q_{i}+Q_{j}\phi(A_{jj})Q_{j}.\eqno(12)$$
For every $A_{ji}\in\mathcal{A}_{ji}$ and $C\in \mathcal{A},$ we have $T_{ji}=\phi^{-1}(A_{ji})\in\mathcal{A}_{ji}.$ Therefore
$$0=\phi([A_{ji},A_{jj}]_{*}\bullet C)=[T_{ji},\phi(A_{jj})]_{*}\bullet\phi(C).$$
Using the surjectivity of $\phi,$ the above equation indicates $[T_{ji},\phi(A_{jj})]_{*}\in\mbox{i}\mathbb{R}I.$
It follows from  Eq. (12) that
$$T_{ji}\phi(A_{jj})Q_{i}-Q_{i}\phi(A_{jj})T_{ji}^{*}\in\mbox{i}\mathbb{R}I.\eqno(13)$$
By  Remark 2.1, multiplying by $Q_{j}$ and $Q_{i}$ from the left and right respectively in Eq. (13), we can get that
$T_{ji}\phi(A_{jj})Q_{i}=0$ for all $T_{ji}\in\mathcal{B}_{ji}.$ By the primeness of $\mathcal{B},$ we obtain that $Q_{i}\phi(A_{jj})Q_{i}=0,$
 thus $\phi(\mathcal{A}_{jj})\subseteq\mathcal{B}_{jj}.$
 Applying the same argument to $\phi^{-1},$ we can obtain ${\mathcal B}_{jj}\subseteq\phi({\mathcal A}_{jj}).$ Consequently, $\phi({\mathcal A}_{jj})={\mathcal B}_{jj}.$

{\bf Lemma 2.11} {\it  $\phi(AB)=\phi(A)\phi(B)$ for all $A,B\in\mathcal{A}.$}

 {\bf Proof} It follows from  Remark 2.1 and Lemma 2.10 that $\phi([P_{i},A_{ij}]_{*}\bullet B_{ji})=[\phi(P_{i}),\phi(A_{ij})]_{*}\bullet\phi(B_{ji})=[Q_{i},\phi(A_{ij})]_{*}\bullet\phi(B_{ji}).$
Thus $$\phi(A_{ij}B_{ji})=\phi(A_{ij})\phi(B_{ji}).\eqno(14)$$
For  $T_{ji}\in\mathcal{B}_{ji},$ we have $X_{ji}=\phi^{-1}(T_{ji})\in\mathcal{A}_{ji}$ by  Lemma 2.10. Therefore
$$\phi(A_{ii}B_{ij})T_{ji}=\phi(A_{ii}B_{ij}X_{ji})=\phi([A_{ii},B_{ij}]_{*}\bullet X_{ji})=\phi(A_{ii})\phi(B_{ij})T_{ji}.$$
 By the primeness of $\mathcal{B},$ we obtain  $$\phi(A_{ii}B_{ij})=\phi(A_{ii})\phi(B_{ij}).\eqno(15)$$
It follows from Eqs. (14)--(15) that
 $$\phi(A_{ij}B_{jj})T_{ji}=\phi(A_{ij}B_{jj}X_{ji})=\phi(A_{ij})\phi(B_{jj}X_{ji})=\phi(A_{ij})\phi(B_{jj})T_{ji}.$$
 In the same manner, we obtain $$\phi(A_{ij}B_{jj})=\phi(A_{ij})\phi(B_{jj}).\eqno(16)$$
By Eq. (15), we have
$$\phi(A_{jj}B_{jj})T_{ji}=\phi(A_{jj}B_{jj}X_{ji})=\phi(A_{jj})\phi(B_{jj}X_{ji})=\phi(A_{jj})\phi(B_{jj})T_{ji}.$$
Thus $$\phi(A_{jj}B_{jj})=\phi(A_{jj})\phi(B_{jj}).\eqno(17)$$
From Eqs. (14)-(17) and  Lemmas 2.9, 2.10, we obtain $\phi(AB)=\phi(A)\phi(B)$ for all $A,B\in\mathcal{A}.$

{\bf Lemma 2.12} {\it   $\phi$  is a linear $*$-isomorphism, or a conjugate linear $*$-isomorphism, or the negative of a linear $*$-isomorphism, or the negative of a conjugate linear $*$-isomorphism.}

 {\bf Proof} It follows from Lemmas 2.9 and 2.11 that  $\phi$  is a ring isomorphism. By Lemma 2.2, exists $\lambda\in\mathbb{R}\setminus\{0\}$
such that $\phi(I)=\lambda I.$ By the equality $\phi(I^{3})=\phi(I)^{3},$ we concur that $\phi(I)=I$ or $\phi(I)=-I.$
In the rest of this section, we deal with these two cases respectively.

{\bf Case 1} If  $\phi(I)=I,$ then  $\phi$  is either a linear $*$-isomorphism or a conjugate linear $*$-isomorphism.

For every rational number $q,$ we have $\phi(qI)=qI.$ Indeed, since $q$ is rational number, exist two integers $r$ and $s$ such that $q=\frac{r}{s}.$
Since $\phi(I)=I$ and  $\phi$  is additive, we get that
$\phi(qI)=\phi(\frac{r}{s}I)=r\phi(\frac{1}{s}I)=\frac{r}{s}\phi(I)=qI.$

Let $A$ be a positive element in $\mathcal{A}.$ Then $A=B^{2}$ for some self-adjoint element $B\in\mathcal{A}.$ It follows from Lemma 2.11 that $\phi(A)=\phi(B)^{2}.$
By Lemma 2.2, we get that $\phi(B)$ is self-adjoint. So $\phi(A)$ is positive. This shows that  $\phi$ preserves positive elements.

Let $\lambda\in\mathbb{R}.$ Choose sequence $\{a_{n}\}$ and $\{b_{n}\}$ of rational numbers such that $a_{n}\leq \lambda \leq b_{n}$ for all $n$ and
$\lim\limits_{n\rightarrow\infty}a_{n}=\lim\limits_{n\rightarrow\infty}b_{n}=\lambda.$ It follows from
$a_{n}I\leq\lambda I\leq b_{n}I$
that
$a_{n}I\leq\phi(\lambda I)\leq b_{n}I.$
Taking the limit, we obtain that $\phi(\lambda I)=\lambda I.$ Hence for all $A\in \mathcal{A},$ we have
$\phi(\lambda A)=\phi((\lambda I)A)=\phi(\lambda I)\phi(A)=\lambda\phi(A).$
Thus $\phi$  is real linear. For every  $A\in\mathcal{A},$ it follows from
$-\phi(A)=\phi(\mbox{i}^{2}A)=\phi(\mbox{i}I)^{2}\phi(A)$ that
$\phi(\mbox{i}I)^{2}=-1,$ which implies that $\phi(\mbox{i}I)=\mbox{i}I$ or $\phi(\mbox{i}I)=-\mbox{i}I.$ By Lemma 2.11, we obtain that
$\phi(\mbox{i}A)=\mbox{i}\phi(A)$ or $\phi(\mbox{i}A)=-\mbox{i}\phi(A)$ for all $A\in \mathcal{A}.$

For all  $A\in \mathcal{A}, A=A_{1}+\mbox{i}A_{2},$ where $A_{1}=\frac{A+A^{*}}{2}$ and $A_{2}=\frac{A-A^{*}}{2\mbox{i}}$ are self-adjoint elements. By Lemmas 2.2 and 2.9,
if $\phi(\mbox{i}A)=\mbox{i}\phi(A),$ then
\begin{eqnarray*}
&&\phi(A^{*})=\phi(A_{1}-\mbox{i}A_{2})=\phi(A_{1})-\phi(\mbox{i}A_{2})=\phi(A_{1})-\mbox{i}\phi(A_{2})\\
&=&\phi(A_{1})^{*}-\mbox{i}\phi(A_{2})^{*}=\phi(A_{1})^{*}+(\mbox{i}\phi(A_{2}))^{*}=\phi(A)^{*}
\end{eqnarray*}
Similarly, if $\phi(\mbox{i}A)=-\mbox{i}\phi(A),$ we also obtain $\phi(A^{*})=\phi(A)^{*}.$
Then  $\phi$  is either a linear $*$-isomorphism or a conjugate linear $*$-isomorphism.

{\bf Case 2} If  $\phi(I)=-I,$ then  $\phi$  is either the negative of  a linear $*$-isomorphism or  the negative of a conjugate linear $*$-isomorphism.

Consider that the map $\psi:\mathcal{A}\rightarrow\mathcal{B}$ defined by $\psi(A)=-\phi(A)$ for all  $A\in \mathcal{A}.$
It is easy to see that $\psi$ satisfies  $\psi([[A,B]_{*},C])=[[\psi(A),\psi(B)]_{*},\psi(C)]$ for all $A,B,C\in\mathcal A$ and $\psi(I)=I.$
Then the arguments for Case 1 ensure that $\psi$  is either a linear $*$-isomorphism or a conjugate linear $*$-isomorphism. So $\phi$ is either the negative of  a linear $*$-isomorphism or  the negative of a conjugate linear $*$-isomorphism.

Combining Cases 1--2, the proof of Theorem 2.1 is finished.

\section{ The mixed Jordan triple $\eta$-$*$-product preserving maps}

{\bf Theorem 3.1}   {\it Let $\mathcal{A},\mathcal{B}$ be two  factor von Neumann algebras and let $\eta\in\mathbb{C}\setminus\{0,1\}.$ Suppose that $\phi$ is a bijective map from
 $\mathcal{A}$ to $\mathcal{B}$ with $\phi([A,B]_{*}^{\eta}\diamond_{\eta} C)=[\phi(A),\phi(B)]_{*}^{\eta}\diamond_{\eta}\phi(C)$ for all $A,B,C\in\mathcal A.$ Then $\phi$  is additive. }

{\bf Theorem 3.2}   {\it Let $\mathcal{A},\mathcal{B}$ be two  factor von Neumann algebras,  let $\eta\in\mathbb{C}\setminus\{0,1\}$ and let
 $\phi:\mathcal{A}\rightarrow\mathcal{B}$ is a bijective map, satisfying  $\phi(I)=I$ and preserving the  mixed Jordan triple $\eta$-$*$-product,
   Then $\phi$ is  either a linear $*$-isomorphism or a conjugate linear $*$-isomorphism.}

{\bf Proof of Theorem 3.1} In the following, we will complete the proof by proving several claims.

{\bf Claim 1} $\phi(0)=0.$

 Since $\phi$ is surjective, there exists $A\in\mathcal{A}$ such that $\phi(A)=0.$ Then we obtain $\phi(0)=\phi([[0,A]_{*}^{\eta}\diamond_{\eta}A)=[[\phi(0),\phi(A)]_{*}^{\eta}\diamond_{\eta}\phi(A)=0.$

{\bf Claim 2}  $\phi(A_{11}+A_{22})=\phi(A_{11})+\phi(A_{22})$ for all $A_{11}\in \mathcal{A}_{11}$ and $A_{22}\in \mathcal{A}_{22}.$

  Let $X=\sum_{i,i=1}^{2}X_{ij}\in\mathcal{A}$  such that $\phi(X)=\phi(A_{11})+\phi(A_{22}).$ For any
  $\lambda\in\mathbb{C},[I,\frac{\lambda P_{1}}{1-\eta}]_{*}^{\eta}\diamond_{\eta}A_{22}=0.$ By applying Eq. (3) and Claim 1, we obtain
 $\phi([I,\frac{\lambda P_{1}}{1-\eta}]_{*}^{\eta}\diamond_{\eta}X)=\phi([I,\frac{\lambda P_{1}}{1-\eta}]_{*}^{\eta}\diamond_{\eta}A_{11}).$
  By the injectivity of $\phi,$ we get that $[I,\frac{\lambda P_{1}}{1-\eta}]_{*}^{\eta}\diamond_{\eta}X=[I,\frac{\lambda P_{1}}{1-\eta}]_{*}^{\eta}\diamond_{\eta}A_{11},$
i.e., $$(\lambda+\overline{\lambda}\eta)X_{11}+\lambda X_{12}+\overline{\lambda}\eta X_{21}=(\lambda+\overline{\lambda}\eta)A_{11}.$$
Assume that $\lambda\neq 0$ and $\lambda+\overline{\lambda}\eta\neq 0,$ we have $X_{11}=A_{11},X_{12}=X_{21}=0.$  In the same manner, we obtain  $X_{22}=A_{22}.$

 {\bf Claim 3}  $\phi(A_{12}+A_{21})=\phi(A_{12})+\phi(A_{21})$ for all $A_{12}\in \mathcal{A}_{12}$ and $A_{21}\in \mathcal{A}_{21}.$

  Let $X=\sum_{i,i=1}^{2}X_{ij}\in\mathcal{A}$  such that $\phi(X)=\phi(A_{12})+\phi(A_{21}).$ For any
  $\lambda\in\mathbb{C},$ since $[I,\frac{\lambda P_{1}-\frac{\overline{\lambda}}{\overline{\eta}}P_{2}}{1-\eta}]_{*}^{\eta}\diamond_{\eta}A_{12}=0,$ Applying Eq. (3) and Claim 1 again, we obtain
 $\phi([I,\frac{\lambda P_{1}-\frac{\overline{\lambda}}{\overline{\eta}}P_{2}}{1-\eta}]_{*}^{\eta}\diamond_{\eta}X)=\phi([I,\frac{\lambda P_{1}-\frac{\overline{\lambda}}{\overline{\eta}}P_{2}}{1-\eta}]_{*}^{\eta}\diamond_{\eta}A_{21}).$
 The injectivity of $\phi$ implies that $[I,\frac{\lambda P_{1}-\frac{\overline{\lambda}}{\overline{\eta}}P_{2}}{1-\eta}]_{*}^{\eta}\diamond_{\eta}X=[I,\frac{\lambda P_{1}-\frac{\overline{\lambda}}{\overline{\eta}}P_{2}}{1-\eta}]_{*}^{\eta}\diamond_{\eta}A_{21},$
i.e., $$(\lambda+\overline{\lambda}\eta)X_{11}-(\lambda+\frac{\overline{\lambda}}{\overline{\eta}})X_{22}+(\overline{\lambda}\eta-\frac{\overline{\lambda}}{\overline{\eta}})X_{21}=(\overline{\lambda}\eta-\frac{\overline{\lambda}}{\overline{\eta}})A_{21}$$
for all $\lambda\in\mathbb{C}.$
Thus we get $X_{11}=X_{22}=0.$

Since $[A_{12},\lambda P_{1}]_{*}^{\eta}\diamond_{\eta}I=0.$    It follows from Eq. (1) that $\phi([[X,\lambda P_{1}]_{*}^{\eta}\diamond_{\eta}I)=\phi([[A_{21},\lambda P_{1}]_{*}^{\eta}\diamond_{\eta}I).$ Thus we obtain $(\lambda-\overline{\lambda}|\eta|^{2})X_{21}+\eta(\overline{\lambda}-\lambda)X_{21}^{*}=(\lambda-\overline{\lambda}|\eta|^{2})A_{21}+\eta(\overline{\lambda}-\lambda)A_{21}^{*}$
for all $\lambda\in\mathbb{C},$  which indicates that $X_{21}=A_{21}.$  In the same manner, we obtain  $X_{12}=A_{12}.$

 {\bf Claim 4}   Let $i,j,k\in\{1,2\}$ with $i\neq j.$ Then
 $\phi(A_{kk}+A_{ij})=\phi(A_{kk})+\phi(A_{ij})$ for all $A_{kk}\in \mathcal{A}_{kk}$ and $A_{ij}\in \mathcal{A}_{ij}.$

 We only prove the case  $i=k=1$ and $j=2,$ the proof of the other cases is similar. Now assume that $X\in\mathcal{A}$ satisfies  $\phi(X)=\phi(A_{11})+\phi(A_{12}).$
 For any
  $\lambda\in\mathbb{C},$ since $[I,\frac{\lambda P_{2}}{1-\eta}]_{*}^{\eta}\diamond_{\eta}A_{11}=0,$ by applying Eq. (3) and Claim 1 again, we obtain
 $\phi([I,\frac{\lambda P_{2}}{1-\eta}]_{*}^{\eta}\diamond_{\eta}X)=\phi([I,\frac{\lambda P_{2}}{1-\eta}]_{*}^{\eta}\diamond_{\eta}A_{12})$ for any $\lambda\in\mathbb{C}.$
 Thus we get  $X_{12}=A_{12},X_{21}=X_{22}=0.$

 Since  $[I,\frac{\lambda P_{1}-\frac{\overline{\lambda}}{\overline{\eta}}P_{2}}{1-\eta}]_{*}^{\eta}\diamond_{\eta}A_{12}=0,$ we have
  $$\phi([I,\frac{\lambda P_{1}-\frac{\overline{\lambda}}{\overline{\eta}}P_{2}}{1-\eta}]_{*}^{\eta}\diamond_{\eta}X)=\phi([I,\frac{\lambda P_{1}-\frac{\overline{\lambda}}{\overline{\eta}}P_{2}}{1-\eta}]_{*}^{\eta}\diamond_{\eta}A_{11})$$  for any $\lambda\in\mathbb{C},$
   which indicates that $X_{11}=A_{11}.$

 {\bf Claim 5}  $\phi(A_{11}+A_{12}+A_{21})=\phi(A_{11})+\phi(A_{12})+\phi(A_{21})$ and $\phi(A_{12}+A_{21}+A_{22})=\phi(A_{12})+\phi(A_{21})+\phi(A_{22})$
 for all $A_{11}\in \mathcal{A}_{11},A_{12}\in \mathcal{A}_{12},A_{21}\in \mathcal{A}_{21}$ and $A_{22}\in \mathcal{A}_{22}.$

 Choose  $X=\sum_{i,i=1}^{2}X_{ij}\in\mathcal{A}$  such that $\phi(X)=\phi(A_{11})+\phi(A_{12})+\phi(A_{21}).$ For any
  $\lambda\in\mathbb{C},$ it follows from Claim 4 and  Eq. (3) that
\begin{eqnarray*}
&&\phi((\lambda+\eta\overline{\lambda})X_{22}+\lambda X_{21}+\overline{\lambda}\eta X_{12})\\
&=&\phi([I,\frac{\lambda P_{2}}{1-\eta}]_{*}^{\eta}\diamond_{\eta}X)\\
&=&\phi([I,\frac{\lambda P_{2}}{1-\eta}]_{*}^{\eta}\diamond_{\eta}A_{11})+\phi([I,\frac{\lambda P_{2}}{1-\eta}]_{*}^{\eta}\diamond_{\eta}A_{12})+\phi([I,\frac{\lambda P_{2}}{1-\eta}]_{*}^{\eta}\diamond_{\eta}A_{21})\\
&=&\phi(\overline{\lambda}\eta A_{12})+\phi(\lambda A_{21})=\phi(\overline{\lambda}\eta A_{12}+\lambda A_{21}),
\end{eqnarray*}
 which indicates that $X_{12}=A_{12},X_{21}=A_{21}$  and $X_{22}=0.$  Thus we get $X=X_{11}+A_{12}+A_{21}.$

 Since  $[I,\frac{-\frac{\overline{\lambda}}{\overline{\eta}}P_{1}+\lambda P_{2}}{1-\eta}]_{*}^{\eta}\diamond_{\eta}A_{21}=0,$  by applying Eq. (3) and Claim 4 again, we obtain
 \begin{eqnarray*}
&&\phi((-\frac{\overline{\lambda}}{\overline{\eta}}-\lambda)X_{11}+(-\frac{\overline{\lambda}}{\overline{\eta}}+\eta\overline{\lambda})X_{12})\\
&=&\phi([I,\frac{-\frac{\overline{\lambda}}{\overline{\eta}}P_{1}+\lambda P_{2}}{1-\eta}]_{*}^{\eta}\diamond_{\eta}X)\\
&=&\phi([I,\frac{-\frac{\overline{\lambda}}{\overline{\eta}}P_{1}+\lambda P_{2}}{1-\eta}]_{*}^{\eta}\diamond_{\eta}A_{11})+\phi([I,\frac{-\frac{\overline{\lambda}}{\overline{\eta}}P_{1}+\lambda P_{2}}{1-\eta}]_{*}^{\eta}\diamond_{\eta}A_{12})\\
&=&\phi((-\frac{\overline{\lambda}}{\overline{\eta}}-\lambda)A_{11})+\phi((-\frac{\overline{\lambda}}{\overline{\eta}}+\eta\overline{\lambda})A_{12})\\
&&=\phi((-\frac{\overline{\lambda}}{\overline{\eta}}-\lambda)A_{11}+(-\frac{\overline{\lambda}}{\overline{\eta}}+\eta\overline{\lambda})A_{12}).
\end{eqnarray*}
This indicates that $X_{11}=A_{11}.$
In the second case, we can similarly prove that the conclusion is valid.

{\bf Claim 6}    Let $i,j\in\{1,2\}$ with $i\neq j.$ Then
 $\phi(A_{ij}+B_{ij})=\phi(A_{ij})+\phi(B_{ij})$ for all $A_{ij},B_{ij}\in \mathcal{A}_{ij}.$

Since $[[I,\frac{P_{i}+A_{ij}}{1-\eta}]_{*}^{\eta},P_{j}+B_{ij}]^{\eta}=A_{ij}+B_{ij}+\eta A_{ij}^{*}+\eta B_{ij}A_{ij}^{*},$   it follows from Claims 5, 4 and 3 that
\begin{eqnarray*}
&&\phi(A_{ij}+B_{ij})+\phi(\eta A_{ij}^{*})+\phi(\eta B_{ij}A_{ij}^{*})=\phi([I,\frac{P_{i}+A_{ij}}{1-\eta}]_{*}^{\eta}\diamond_{\eta}(P_{j}+B_{ij}))\\
&=&[\phi(I),\phi(\frac{P_{i}+A_{ij}}{1-\eta})]_{*}^{\eta}\diamond_{\eta}\phi(P_{j}+B_{ij})\\
&=&[\phi(I),\phi(\frac{P_{i}}{1-\eta})+\phi(\frac{A_{ij}}{1-\eta})]_{*}^{\eta}\diamond_{\eta}(\phi(P_{j})+\phi(B_{ij})))\\
&=&[\phi(I),\phi(\frac{P_{i}}{1-\eta})]_{*}^{\eta}\diamond_{\eta}\phi(P_{j})+[\phi(I),\phi(\frac{P_{i}}{1-\eta})]_{*}^{\eta}\diamond_{\eta}\phi(B_{ij})\\
&&+[\phi(I),\phi(\frac{A_{ij}}{1-\eta})]_{*}^{\eta}\diamond_{\eta}\phi(P_{j})+[\phi(I),\phi(\frac{A_{ij}}{1-\eta})]_{*}^{\eta}\diamond_{\eta}\phi(B_{ij})\\
&=&\phi(B_{ij})+\phi(A_{ij}+\eta A_{ij}^{*})+\phi(\eta B_{ij}A_{ij}^{*})\\
&=&\phi(B_{ij})+\phi(A_{ij})+\phi(\eta A_{ij}^{*})+\phi(\eta B_{ij}A_{ij}^{*}).
\end{eqnarray*}
 Then  $\phi(A_{ij}+B_{ij})=\phi(A_{ij})+\phi(B_{ij}).$

 {\bf Claim 7}  $\phi(T_{12}A_{21}+T_{12}B_{21}+\eta A_{12}T_{12}^{*}+\eta B_{12}T_{12}^{*})=\phi(T_{12}A_{21})+\phi(T_{12}B_{21})+\phi(\eta A_{12}T_{12}^{*})+\phi(\eta B_{12}T_{12}^{*})$
 for all $T_{12},A_{12},B_{12}\in \mathcal{A}_{12}$ and  $A_{21},B_{21}\in \mathcal{A}_{21}.$

 It follows from Claims 3 and 6 that
\begin{eqnarray*}
&&\phi(T_{12}A_{21}+T_{12}B_{21}+\eta A_{12}T_{12}^{*}+\eta B_{12}T_{12}^{*})\\
&=&\phi([[I,\frac{T_{12}}{1-\eta}]_{*}^{\eta}\diamond_{\eta}(A_{21}+B_{21}+A_{12}+B_{12}))\\
&=&[\phi(I),\phi(\frac{T_{12}}{1-\eta})]_{*}^{\eta}\diamond_{\eta}(\phi(A_{21})+\phi(B_{21})+\phi(A_{12})+\phi(B_{12}))\\
&=&\phi([I,\frac{T_{12}}{1-\eta}]_{*}^{\eta}\diamond_{\eta}A_{21})+\phi([I,\frac{T_{12}}{1-\eta}]_{*}^{\eta}\diamond_{\eta}B_{21})\\
&&+\phi([I,\frac{T_{12}}{1-\eta}]_{*}^{\eta}\diamond_{\eta}A_{12})+\phi([I,\frac{T_{12}}{1-\eta}]_{*}^{\eta}\diamond_{\eta}B_{12})\\
&=&\phi(T_{12}A_{21})+\phi(T_{12}B_{21})+\phi(\eta A_{12}T_{12}^{*})+\phi(\eta B_{12}T_{12}^{*}).
\end{eqnarray*}

 {\bf Claim 8}  $\phi([I,\frac{T_{12}}{1-\eta}]_{*}^{\eta}\diamond_{\eta}(A+B))=\phi([I,\frac{T_{12}}{1-\eta}]_{*}^{\eta}\diamond_{\eta}A)+\phi([I,\frac{T_{12}}{1-\eta}]_{*}^{\eta}\diamond_{\eta}B)$
 for all $T_{12}\in \mathcal{A}_{12}$ and  $A,B\in \mathcal{A}.$

 Let $A=\sum_{i,i=1}^{2}A_{ij}$  and $B=\sum_{i,i=1}^{2}B_{ij},$   it follows from Claims 5, 6 and  7 that
\begin{eqnarray*}
&&\phi([I,\frac{T_{12}}{1-\eta}]_{*}^{\eta}\diamond_{\eta}(A+B))\\
&=&\phi(T_{12}A_{22}+T_{12}B_{22}+T_{12}A_{21}+T_{12}B_{21}\\
&&+\eta A_{12}T_{12}^{*}+\eta B_{12}T_{12}^{*}+\eta A_{22}T_{12}^{*}+\eta B_{22}T_{12}^{*})\\
&=&\phi(T_{12}A_{22}+T_{12}B_{22})+\phi(T_{12}A_{21}+T_{12}B_{21}+\eta A_{12}T_{12}^{*}+\eta B_{12}T_{12}^{*})\\
&&+\phi(\eta A_{22}T_{12}^{*}+\eta B_{22}T_{12}^{*})\\
&=&\phi(T_{12}A_{22})+\phi(T_{12}B_{22})+\phi(T_{12}A_{21})+\phi(T_{12}B_{21})\\
&&+\phi(\eta A_{12}T_{12}^{*})+\phi(\eta B_{12}T_{12}^{*})+\phi(\eta A_{22}T_{12}^{*})+\phi(\eta B_{22}T_{12}^{*})\\
&=&\phi(T_{12}A_{22})+\phi(T_{12}A_{21}+\eta A_{12}T_{12}^{*})+\phi(\eta A_{22}T_{12}^{*})\\
&&+\phi(T_{12}B_{22})+\phi(T_{12}B_{21}+\eta B_{12}T_{12}^{*})+\phi\eta B_{22}T_{12}^{*})\\
&=&\phi(T_{12}A_{22}+T_{12}A_{21}+\eta A_{12}T_{12}^{*}+\eta A_{22}T_{12}^{*})\\
&&+\phi(T_{12}B_{22}+T_{12}B_{21}+\eta B_{12}T_{12}^{*}+\eta B_{22}T_{12}^{*})\\
&=&\phi([I,\frac{T_{12}}{1-\eta}]_{*}^{\eta}\diamond_{\eta}A)+\phi([I,\frac{T_{12}}{1-\eta}]_{*}^{\eta}\diamond_{\eta}B).
\end{eqnarray*}

{\bf  Claim 9}  $\phi$ is additive.

Let $A$ and $B$ be in  $\mathcal{A}.$ To show that $\phi(A+B)=\phi(A)+\phi(B),$   we choose $X\in\mathcal{A}$ such that $\phi(X)=\phi(A)+\phi(B).$
For any  $T_{12}\in \mathcal{A}_{12},$  it follows from Claim 8 and  Eq. (3) that
\begin{eqnarray*}
\phi([I,\frac{T_{12}}{1-\eta}]_{*}^{\eta}\diamond_{\eta}X)&=&\phi([I,\frac{T_{12}}{1-\eta}]_{*}^{\eta}\diamond_{\eta}A)+\phi([I,\frac{T_{12}}{1-\eta}]_{*}^{\eta}\diamond_{\eta}B)\\
&=&\phi([I,\frac{T_{12}}{1-\eta}]_{*}^{\eta}\diamond_{\eta}(A+B)),
\end{eqnarray*}
 which indicates that $$[I,\frac{T_{12}}{1-\eta}]_{*}^{\eta}\diamond_{\eta}X=[I,\frac{T_{12}}{1-\eta}]_{*}^{\eta}\diamond_{\eta}(A+B).$$
Thus we get $T_{12}(X-A-B)-\eta(X-A-B)T_{12}^{*}=0.$ By $\mbox{i}T_{12}$ replacing $T_{12}$ in the above equation, we get $T_{12}(X-A-B)+\eta(X-A-B)T_{12}^{*}=0$
and hence $T_{12}(X-A-B)=0$ for all  $T_{12}\in \mathcal{A}_{12},$ this indicates that $X-A-B=0.$
  We concur that $\phi$ is additive, and the proof of Theorem 3.1 is completed.

{\bf Proof of Theorem 3.2} In the rest of this section, we will deal with two cases.

{\bf Case 1}  $|\eta|=1.$

For all $A,B,C\in\mathcal{A},$  it follows from
$$\phi([A,B]_{*}^{\eta}\diamond_{\eta}C)=[\phi(A),\phi(B)]_{*}^{\eta}\diamond_{\eta}\phi(C)$$
that
\begin{eqnarray*}
&&\phi(ABC-\eta BA^{*}C+\eta CB^{*}A^{*}-|\eta|^{2}CAB^{*})\\
&=&\phi(A)\phi(B)\phi(C)-\eta \phi(B)\phi(A)^{*}\phi(C)\\
&&+\eta\phi(C)\phi(B)^{*}\phi(A)^{*}-|\eta|^{2}\phi(C)\phi(A)\phi(B)^{*}. \ \ \ \ \ \ \ \ \ \ \  \ \ \ \ \ \ \ \ \ \ \  \ \ \ \ \ \ \ \  \ (18)
\end{eqnarray*}
Replacing $\eta$ with $-\eta$ in Eq. (18), we obtain
\begin{eqnarray*}
&&\phi(ABC+\eta BA^{*}C-\eta CB^{*}A^{*}-|\eta|^{2}CAB^{*})\\
&=&\phi(A)\phi(B)\phi(C)+\eta\phi(B)\phi(A)^{*}\phi(C)\\
&&-\eta\phi(C)\phi(B)^{*}\phi(A)^{*}-|\eta|^{2}\phi(C)\phi(A)\phi(B)^{*}.  \ \ \ \ \ \ \ \ \ \ \  \ \ \ \ \ \ \ \ \ \ \  \ \ \ \ \ \ \ \ \ (19)
\end{eqnarray*}
It follows from Eqs. (18)--(19) and  Theorem 3.1 that
$$\phi(ABC-CAB^{*})=\phi(A)\phi(B)\phi(C)-\phi(C)\phi(A)\phi(B)^{*}.\eqno(20)$$
Taking $A=I$ in  Eq. (20), since $\phi(I)=I,$ we obtain
$$\phi(BC-CB^{*})=\phi(B)\phi(C)-\phi(C)\phi(B)^{*}$$  for all $B,C\in\mathcal{A}.$
Based on the result of [9],  $\phi$ is a $*$-ring isomorphism.

{\bf Case 2}  $|\eta|\neq 1.$

Take $\alpha=\frac{1}{1-|\eta|^{2}}.$ Then we obtain $\alpha(1-|\eta|^{2})=1$ and $\alpha\neq 0.$

{\bf Clam 2.1} $\phi(\alpha I)=\alpha I$ and $\phi$ preserves self-adjoint elements in both directions.

Since $I=\phi([\alpha I,I]_{*}^{\eta}\diamond_{\eta}I)=[\phi(\alpha I),I]_{*}^{\eta}\diamond_{\eta}I=(1-|\eta|^{2})\phi(\alpha I)=\frac{1}{\alpha}\phi(\alpha I),$
we have $\phi(\alpha I)=\alpha I.$

Let $A\in\mathcal{A}$ such that $A=A^{*}.$ Then
\begin{eqnarray*}
&&\phi(A)=\phi([\alpha I,A]_{*}^{\eta}\diamond_{\eta}I)=[\alpha I,\phi(A)]_{*}^{\eta}\diamond_{\eta}I)\\
&=&(\alpha-\overline{\alpha}\eta)\phi(A)+(\overline{\alpha}\eta-\alpha|\eta|^{2})\phi(A)^{*}.
\end{eqnarray*}
Therefore $$\phi(A)^{*}=\frac{1-\alpha+\overline{\alpha}\eta}{\overline{\alpha}\eta-\alpha|\eta|^{2}}\phi(A)=\phi(A),$$
which indicates the sufficiency. The necessity can be obtained by considering $\phi^{-1}.$

{\bf Clam 2.2} $\frac{1}{\alpha}\phi(\alpha P)$ is a projection in $\mathcal{B}$ if and only if $P$ is  a projection in $\mathcal{A}.$

Since $(\alpha P)^{*}=\alpha P,$ it follows from Claim 2.1 that  $\phi(\alpha P)^{*}=\phi(\alpha P).$ Thus $\frac{1}{\alpha}\phi(\alpha P)$ is
self-adjoint. Moreover,
\begin{eqnarray*}
&&\phi(\alpha P)=\phi([\alpha P,I]_{*}^{\eta}\diamond_{\eta}\alpha P)\\
&=&[\phi(\alpha P),I]_{*}^{\eta}\diamond_{\eta}\phi(\alpha P)=(1-|\eta|^{2})\phi(\alpha P)^{2}=\frac{1}{\alpha}\phi(\alpha P)^{2}
\end{eqnarray*}
and $(\frac{1}{\alpha}\phi(\alpha P))^{2}=\frac{1}{\alpha}\phi(\alpha P).$ So  $\frac{1}{\alpha}\phi(\alpha P)$ is a projection, which shows the sufficiency.
The  necessity can be proved by considering $\phi^{-1}.$

{\bf Clam 2.3}  Let $i,j\in\{1,2\}$ with $i\neq j.$ Then $\phi(\mathcal{A}_{ij})=\mathcal{B}_{ij}.$

Choose  a projection $Q_{i}\in\mathcal{B},i=1,2,$ by Clam 2.2, we have $P_{i}=\frac{1}{\alpha}\phi^{-1}(\alpha Q_{i})$ is projection in $\mathcal{A}.$
It is easy to see that $Q_{i}=\frac{1}{\alpha}\phi(\alpha P_{i}).$

For any $A_{ij}\in\mathcal{A}_{ij},$ since
\begin{eqnarray*}
&&\phi(\alpha(1-\eta)A_{ij})=\phi([I,\alpha P_{i}]_{*}^{\eta}\diamond_{\eta}A_{ij})\\
&=&[I,\alpha Q_{i}]_{*}^{\eta}\diamond_{\eta}\phi(A_{ij})=\alpha(1-\eta)Q_{i}\phi(A_{ij})+\overline{\alpha}\eta(1-\overline{\eta})\phi(A_{ij})Q_{i},
\end{eqnarray*}
we obtain $Q_{j}\phi(\alpha(1-\eta)A_{ij})Q_{j}=0.$  In the same manner, we obtain  $Q_{i}\phi(\overline{\alpha}\eta(\overline{\eta}-1)A_{ij})Q_{i}=0.$
Since $A_{ij}$ is arbitrary, we obtain $\phi(A_{ij})=B_{ij}+B_{ji}$ for some  $B_{ij}\in\mathcal{B}_{ij}$ and $B_{ji}\in\mathcal{B}_{ji}.$
Because
\begin{eqnarray*}
&&0=\phi([I,A_{ij}]_{*}^{\eta}\diamond_{\eta}\alpha P_{i})\\
&=&[I,\phi(A_{ij})]_{*}^{\eta}\diamond_{\eta}\alpha Q_{i})=\alpha(1-\eta)B_{ji}+\alpha\eta(1-\overline{\eta})B_{ji}^{*}
\end{eqnarray*}
and $\alpha(1-\eta)\neq 0,$ we obtain $B_{ji}=0,$ this indicates $\phi(A_{ij})\in\mathcal{B}_{ij}.$ By considering $\phi^{-1},$ we get  $\phi(\mathcal{A}_{ij})=\mathcal{B}_{ij}.$

{\bf Clam 2.4} $\phi(\mathcal{A}_{ii})\subseteq\mathcal{B}_{ii},i=1,2.$

Let $A_{ii}\in\mathcal{A}_{ii}$ and $i\neq j,$ we have
\begin{eqnarray*}
&&0=\phi([I,\alpha P_{j}]_{*}^{\eta}\diamond_{\eta}A_{ii})\\
&=&[I,\alpha Q_{j}]_{*}^{\eta}\diamond_{\eta}\phi(A_{ii})=\alpha(1-\eta)Q_{j}\phi(A_{ii})+\overline{\alpha}\eta(1-\overline{\eta})\phi(A_{ii})Q_{j},
\end{eqnarray*}
which indicates that $Q_{j}\phi(A_{ii})Q_{i}=Q_{i}\phi(A_{ii})Q_{j}=0$ and $\phi(A_{ii})=B_{ii}+B_{jj}$ for some  $B_{ii}\in\mathcal{B}_{ii}$ and $B_{jj}\in\mathcal{B}_{jj}.$

Let $T_{ij}\in\mathcal{B}_{ij}$ and $i\neq j.$ It follows from Claim 2.3 that $\phi^{-1}(T_{ij})\in\mathcal{A}_{ij},$ thus
\begin{eqnarray*}
&&0=\phi([I,\phi^{-1}(T_{ij})]_{*}^{\eta}\diamond_{\eta}A_{ii})\\
&=&[I,T_{ij}]_{*}^{\eta}\diamond_{\eta}\phi(A_{ii})=(1-\eta)T_{ij}B_{jj}+\eta(1-\overline{\eta})B_{jj}T_{ij}^{*},
\end{eqnarray*}
which indicates that $B_{jj}=0.$ So we have $\phi(A_{ii})=B_{ii}\subseteq\mathcal{B}_{ii}.$

{\bf Clam 2.5} $\phi(AB)=\phi(A)\phi(B)$ for all $A,B\in\mathcal{A}.$

To prove  $\phi(AB)=\phi(A)\phi(B),$ we need to consider that $\phi(A_{ij}B_{kl})=\phi(A_{ij})\phi(B_{kl})$ for any $i,j,k,l\in\{1,2\}.$ If $j\neq k,$ it follows from
Claims 2.3 and 2.4 that $\phi(A_{ij}B_{kl})=\phi(A_{ij})\phi(B_{kl})=0,$ then we just need to prove the cases with $j=k.$

It follows from $\phi(B_{ij})\phi(A_{ii})^{*}=0$ that
\begin{eqnarray*}
&&\phi(A_{ii}B_{ij})+\phi(\eta B_{ij}^{*}A_{ii}^{*})=\phi([A_{ii},B_{ij}]_{*}^{\eta}\diamond_{\eta}I)\\
&=&[\phi(A_{ii}),\phi(B_{ij})]_{*}^{\eta}\diamond_{\eta}I)=\phi(A_{ii})\phi(B_{ij})+\eta\phi(B_{ij})^{*}\phi(A_{ii})^{*}.
\end{eqnarray*}
By Claims 2.3 and 2.4, we obtain $\phi(A_{ii}B_{ij})=\phi(A_{ii})\phi(B_{ij}).$ 

For any $T_{ij}\in\mathcal{B}_{ij},i\neq j,$ we have $C_{ij}=\phi^{-1}(T_{ij})\in\mathcal{A}_{ij}$ by  Clam 2.3. So
$$\phi(A_{ii}B_{ii})T_{ij}=\phi(A_{ii}B_{ii}C_{ij})=\phi(A_{ii})\phi(B_{ii}C_{ij})=\phi(A_{ii})\phi(B_{ii})T_{ij}.$$
By the primeness of $\mathcal{B},$ we obtain that $\phi(A_{ii}B_{ii})=\phi(A_{ii})\phi(B_{ii}).$

Since $\phi(B_{ji})\phi(A_{ij})^{*}=0,$ we have
\begin{eqnarray*}
&&\phi(A_{ij}B_{ji})+\phi(\eta B_{ji}^{*}A_{ij}^{*})=\phi([A_{ij},B_{ji}]_{*}^{\eta}\diamond_{\eta}I)\\
&=&[\phi(A_{ij}),\phi(B_{ji})]_{*}^{\eta}\diamond_{\eta}I)=\phi(A_{ij})\phi(B_{ji})+\eta\phi(B_{ji})^{*}\phi(A_{ij})^{*}.
\end{eqnarray*}
which indicates that $\phi(A_{ij}B_{ji})=\phi(A_{ij})\phi(B_{ji}).$

For any $T_{ji}\in\mathcal{B}_{ji},i\neq j,$ we have $C_{ji}=\phi^{-1}(T_{ji})\in\mathcal{A}_{ji}.$  So
$$\phi(A_{ij}B_{jj})T_{ji}=\phi(A_{ij}B_{jj}C_{ji})=\phi(A_{ij})\phi(B_{jj}C_{ji})=\phi(A_{ij})\phi(B_{jj})T_{ji}.$$
This indicates that $\phi(A_{ij}B_{jj})=\phi(A_{ij})\phi(B_{jj}).$

Combining Cases 1--2,  similar to the case 1 of  Theorem 2.1, the proof of Theorem 3.2 is finished.

\end{document}